\documentclass[12pt]{article}
\usepackage[T2A]{fontenc}
\usepackage[cp1251]{inputenc}
\usepackage[english]{babel}
\usepackage{amsmath,amsfonts,amssymb, color}
\newtheorem{Theorem}{Theorem}
\newtheorem{Remark}{Remark}
\newtheorem{Lemma}{Lemma}
\newtheorem{Proposition}{Proposition}
\newtheorem{Corollary}{Corollary}

\title{\bf Spectral analysis 
for some multifractional Gaussian processes}
\author{A.I. Karol\footnote{St.Petersburg State University, Universitetskaya emb. 7-9, St.Petersburg, 199034, Russia; E-mail: andrey.i.karol@gmail.com}\setcounter{footnote}{6}
\ and A.I. Nazarov\footnote{ St.Petersburg Department of Steklov Mathematical Institute of Russian Academy of Science, Fontanka 27, St.Petersburg, 191023, Russia, and St.Petersburg State University, Universitetskaya emb. 7-9, St.Petersburg, 199034, Russia; E-mail: nazarov@pdmi.ras.ru. Partially supported by the Russian Foundation of Basic Research Grant 20-51-12004.}
}

\begin{document}

\maketitle

\hfill {\it Dedicated to the memory of Ya.Yu. Nikitin,}

\hfill {\it our colleague and friend}
 \medskip
 
 \begin{abstract}
  We study the small ball asymptotics problem in $L_2$ for two generalizations of the fractional Brownian motion with variable Hurst parameter. To this end, we perform careful analysis of the singular values asymptotics for associated integral operators.
 \end{abstract}

 \section{Introduction}

The spectral asymptotics for Gaussian processes are intensively investigated in the last two decades, closely connected with the problem of small deviation asymptotics for such processes in the Hilbert norm. Namely, see \cite{Naz09b}, to obtain the {\it logarithmic} $L_2$-small ball asymptotics of a Gaussian process $X$, it is sufficient to know one-term asymptotics of the eigenvalue counting function of its covariance operator. 

In the most elaborated case of the so-called {\it Green Gaussian processes}, i.e. the processes the covariance functions ${\cal G}_X$ of which are the Green functions for the ordinary differential operators (ODO), one can obtain even two-term asymptotics of the eigenvalues with the remainder estimate and thus manage the {\it exact} small deviation asymptotics (up to a constant). This approach was developed in \cite{NaNi04}, \cite{Naz09a}, see more references in \cite{NaNi18}. 

The case of {\it fractional Gaussian processes} is more complicated. In the pioneer paper \cite{Bron} the one-term spectral asymptotics was calculated for the fractional Brownian motion (FBM) $W^H$, i.e. the zero mean-value Gaussian process with covariance function
\begin{equation*}
{\cal G}_{W^H}(x,y)=\frac 12\,\big(x^{2H}+y^{2H}-|x-y|^{2H}\big)
\end{equation*}
(here $H\in(0,1)$ is the so-called {\it Hurst index}, the case $H=\frac 12$ corresponds to the standard Wiener process). 

A more general approach was suggested in \cite{NaNi04a}. This approach is based on the powerful theorems on spectral asymptotics of integral operators \cite{BS70}, see also \cite[Appendix 7]{BS74}, 
and covers many fractional processes. Now the problem of logarithmic
 $L_2$-small ball asymptotics for such processes is also well-studied. We mention also the breakthrough paper \cite{ChiKl} where the two-term asymptotics of the eigenvalues with the remainder estimate was obtained for the FBM in the full range of the Hurst index. Some generalization of the seminal idea of \cite{ChiKl} was given in \cite{Naz20} where the reader also can find an extensive bibliography.
\medskip

In this paper we consider some more sophisticated Gaussian processes. 

The {\it multifractional Brownian Motion} (mBM) was introduced in  \cite{PV} and \cite{BCI} and was investigated in several papers, see, e.g., \cite{Coh}, \cite{ACV}, and \cite{co}. There are some different definitions of mBM equivalent up to a multiplicative deterministic function. We choose the so-called harmonizable representation \cite{BCI}
\begin{equation}\label{mBM}
W^{H(\cdot)}(x)=C_*(H(x)) \int\limits_{-\infty}^{\infty}\frac{e^{ix\xi}-1}{|\xi|^{H(x)+\frac 12}}\,dW(\xi),
\end{equation} 
where $W(\xi)$ is a conventional Wiener process and the functional Hurst parameter $H(x)$ satisfies $0<H(x)<1$. The choice of normalizing factor
\begin{equation*}
C_*(H)= \left( \frac{\Gamma(2H +1)\sin(\pi H)}{2\pi}\right) ^{\frac 12}
\end{equation*}
ensures that $EX^2(1)=1$.
\medskip

A different process of the same structure is the {\it multifractal Brownian Motion} (mfBM) introduced in \cite{RS}, see also \cite{RJ}. It can be constructed as
\begin{equation}\label{mfBM}
X^{H(\cdot)}(x)= \int\limits_0^x K(x,y,H(x))\,dW(y),
\end{equation} 
where
\begin{equation}\label{K}
K(x,y,H)= c_*(H) y^{\frac 12-H}\int\limits_y^x (z-y)^{H-\frac 32}z^{H-\frac 12}\, dz\, {\mathbb I}_{[0,x]}(y),
\end{equation}
while the variable Hurst parameter $H(x)$ satisfies $\frac 12<H(x)<1$. The normalizing factor is defined by the formula
\begin{equation*}
c_*(H)=\left(\frac{H(2H-1)\Gamma(\frac 32-H)}{\Gamma(2-2H)\Gamma(H-\frac 12)}\right) ^{\frac 12}, 
\end{equation*}
where $\Gamma$ is Euler's gamma-function. 

Both processes (\ref{mBM}) and (\ref{mfBM}) obviously have zero mean, their covariance functions were derived in \cite{ACV} and \cite{RJ}, respectively. For $H(x)\equiv H=const$ they both coincide with conventional FBM.
\medskip 

 We derive one-term spectral asymptotics for the processes (\ref{mBM}) and (\ref{mfBM}) under some regularity assumptions on the functional parameter $H(x)$. Then, using the results of \cite{KNN} we obtain the logarithmic $L_2$-small ball asymptotics for these processes. Despite the fact that the behavior of covariances of mBM and mfBM is significantly different, it turns out that under the assumption $H(x)>\frac 12$ these logarithmic asymptotics coincide.
 \medskip
 
 The structure of our paper is as follows. In Section \ref{S:2} we introduce operators associated with processes under consideration and formulate the result concerning the asymptotics of their singular values. This result is proved in Section \ref{S:3}. Section \ref{S:4} is devoted to $L_2$-small ball behavior of mBM and mfBM. Some auxiliary estimates and asymptotics of singular values of compact operators are collected in Appendix.
 \medskip
  
 We use the letter $C$ to denote various positive constants. To indicate that $C$ depends on some parameters, we list them in parentheses: $C(\dots)$.

\section{Operators associated with mBM and mfBM}
\label{S:2}

Here we define integral operators associated with processes (\ref{mBM}) and (\ref{mfBM}), see \cite[\S 3.2]{Lif}:
\begin{equation}\label{T}
{\mathbb{T}}: L_2(\mathbb{R})\to L_2(0,1),\qquad (\mathbb{T}f)(x):= C_*(H(x)) \int\limits_{-\infty}^{\infty}\frac{e^{ix\xi}-1}{|\xi|^{H(x)+\frac 12}}\,f(\xi)\,d\xi,
\end{equation}
\begin{equation}\label{S} 
{\mathbb{S}}: L_2(0,1)\to L_2(0,1),\qquad ({\mathbb{S}}f)(x):= \int\limits_0^x K(x,y,H(x))f(y)dy
\end{equation}
(the function $K$ is defined in (\ref{K})). It is easy to see that the covariance functions
$$
{\cal G}_{W^{H(\cdot)}}(x,y):=\mathbb{E}W^{H(\cdot)}(x)W^{H(\cdot)}(y),\qquad {\cal G}_{X^{H(\cdot)}}(x,y):=\mathbb{E}X^{H(\cdot)}(x)X^{H(\cdot)}(y)
$$
are the kernels of integral operators $\mathbb{T}\mathbb{T}^*$ and 
$\mathbb{S}\mathbb{S}^*$, respectively.
\medskip

In what follows we denote by $\{\lambda_k (\mathbb{K})\}$ the nonincreasing sequence of eigenvalues of a compact selfadjoint positive operator $\mathbb{K}$ in a Hilbert space $\mathcal H$,
enumerated with multiplicities.

Recall that for a compact operator $\mathbb{A} : \mathcal{H}_1 \to \mathcal{H}_2$, we have $\lambda_k (\mathbb{A}^*\mathbb{A})=\lambda_k (\mathbb{A}\mathbb{A}^*)$, and square roots of these eigenvalues $s_k(\mathbb{A}):= (\lambda_k (\mathbb{A}^*\mathbb{A}))^{\frac 12}$ are called {\it the singular values} of
the operator $\mathbb{A}$.

We denote by ${\cal N}(t,\mathbb{A})$, $t>0$, the {\it distribution function} of singular values,
$${\cal N}(t,\mathbb{A}):=\#\{k \mid s_k(\mathbb{A})>t^{-1}\}. 
$$
Notice that the function $t\mapsto {\cal N}(t,\mathbb{A})$ is conceptually inverse to the function $k\mapsto s_k^{-1}(\mathbb{A})$. Thus, if the singular values have moderate speed of decay then the one-term asymptotic of $s_k(\mathbb{A})$ as $k\to\infty$ is uniquely defined by the one-term asymptotic of ${\cal N}(t,\mathbb{A})$ as $t\to\infty$.\medskip

We suppose that the functional parameter $H(x)$ is a H\"older continuous function. Then it turns out that the singular value asymptotics for operators (\ref{T}), (\ref{S}) heavily depend on the behavior of $H(x)$ in a neighborhood of the set where it attains its minimal value. We set 
$$
H_{min}:=\min\limits_{x\in[0,1]} H(x), \qquad
{\bf D}:=\{x \in [0,1]\mid H(x)=H_{min}\}. 
$$
In what follows we use the notation $\mathfrak{m}=H_{min}+\frac 12$.
\medskip

Next, we introduce the {\it regularized distance} to ${\bf D}$, see, e.g., \cite{Lie}, that is a function ${\rm d}(x)$, $x \in [0,1]$, such that ${\rm d}(x) \asymp {\rm dist}(x,{\bf D})$ and
\begin{equation*} 
{\rm d}\in {\cal C}^\infty([0,1]\setminus {\bf D}),\qquad
|{\rm d}^{(n)}(x)| \le C(n) {\rm d}^{1-n}(x), \quad x
\in [0,1]\setminus {\bf D}, \ \ n\in\mathbb{N}.
\end{equation*}

We describe the behavior of $H(x)$ in a neighborhood of ${\bf D}$ by the following assumptions:
\smallskip

{\bf 1}. The function $h(x):=H(x)-H_{min}$ is bounded by a power of the distance ${\rm d}(x)$:
$$h(x) \le C {\rm d}^\kappa(x) \qquad \mbox{for some} \quad \kappa >0.
$$

{\bf 2}. The function $h(x)$ admits an asymptotic representation with the smooth main term in a neighborhood of ${\bf D}$. More precisely, $h(x)=h_0(x)+h_1(x)$, where $h_0, h_1\in {\cal C}^\beta [0,1]$ with $\beta >0$, and
$$h_0 \in {\cal C}^\infty\big([0,1] \setminus {\bf D}\big), \quad |h^{(n)}_0 (x)|\le C(n)
h_0(x){\rm d}^{-n}(x), \quad x \in [0,1] \setminus {\bf D},$$
\begin{equation*}
h_1(x)=O\big(h^{1+\tau}_0(x)\big) \quad\mbox{as} \quad h_0(x) \to 0\qquad \mbox{for some} \quad \tau >0.
\end{equation*}

{\bf 3}. The measure of the small values set for the function $h$ is a regularly varying function of the level:
\begin{equation}\label{H0} 
{\rm meas}\,\{x\in [0,1] \mid 0<h(x)<s\} = s^\sigma \varphi(s^{-1}),\quad 0<s<s_0,
\end{equation}
where $\sigma, s_0>0$ and $\varphi$ is a {\it slowly varying function} (SVF), see \cite{Sen}.
\smallskip

\begin{Theorem}\label{T:1} 
Let the functional parameter $H(x)$ satisfy the assumptions {\bf 1}--{\bf 3}. Assume in addition that
\begin{equation}\label{tau}
\tau>
\begin{cases}
 \max\{0, \frac 12\,(1-\sigma)\}, & if \quad {\rm meas}\, {\bf D}>0;\\
 \sigma \mathfrak{m}+\frac 12\,(1-\sigma), & if \quad {\rm meas}\, {\bf D}=0.
\end{cases}
\end{equation}
Then, as $t \to \infty$,
\begin{equation}\label{TT1}
{\cal N}(t,\mathbb{T})=\frac 1{\pi}\,\Big((2\pi)^{\frac 12}C_* (H_{min})\Big)^{\frac 1{\mathfrak{m}}} \int\limits_0^1 t^{\frac 1{H(x)+\frac 12\vphantom{1^{1^1}}}}\,dx\,\big(1+O(\log^{-\nu}t)\big);  
\end{equation}
\begin{equation}\label{SS1}
{\cal N}(t,\mathbb{S})=
\frac 1{\pi}\,\Big(c_* (H_{min})\Gamma\big(H_{min}-\frac 12\big)\Big)^{\frac 1{\mathfrak{m}}}
\int\limits_0^1 t^{\frac 1{H(x)+\frac 12\vphantom{1^{1^1}}}}\,dx\,\big(1+O(\log^{-\nu}t)\big),  
\end{equation}
Here $\nu =\nu (H_{min},\sigma,\tau)>0$. If ${\rm meas}\, {\bf D}>0$ then $\nu$ is arbitrary exponent less than $\frac{\tau -\frac 12\,(1-\sigma)}{{\mathfrak{m}}+1}$.
\end{Theorem}

\begin{Remark}
\label{R:1}
 \begin{enumerate}
  \item We recall that the definition of the operator
 $\mathbb{T}$ admits $0<H(x)<1$ while the definition of the operator $\mathbb{S}$ admits only $\frac 12<H_{min}<1$. However, for $H_{min}>\frac 12$ asymptotic formulae (\ref{TT1}) and (\ref{SS1}) coincide since
 $$
 (2\pi)^{\frac 12}C_* (H)=c_* (H)\Gamma\big(H-\frac 12\big)=\left( \Gamma(2H +1)\sin(\pi H)\right) ^{\frac 12}.
 $$
 
  \item If $h_1(x) \equiv 0$, then formulae (\ref{TT1}) and (\ref{SS1}) are valid with $O(t^{-r})$ remainder estimate for some $r>0$.

  \item Applying Laplace method for the integral in (\ref{TT1}) and (\ref{SS1}), we obtain, as $t\to\infty$,
\begin{equation}\label{Laplace}
\int\limits_0^1 t^{\frac 1{H(x)+\frac 12\vphantom{1^{1^1}}}}\,dx =  t^{\frac 1{\mathfrak{m}}}
\big({\rm meas}\, {\bf D} + \mathfrak{m}^{2\sigma} \Gamma(\sigma+1)(\log t)^{-\sigma}  \, \varphi (\log t)(1+o(1))\big).
\end{equation}
Notice that if ${\rm meas}\, {\bf D} >0$ then the main term of the asymptotics is purely power. If in addition $\sigma<\nu$ then we have even two-term asymptotics. In the case $\nu \le \sigma$ we obtain only one-term power asymptotics with logarithmic remainder term.

\item We stress that the asymptotic (\ref{Laplace}) does not change if we replace the function $\varphi$ with an equivalent SVF.
 \end{enumerate}
\end{Remark}

\section{Proof of Theorem 1}
\label{S:3}

The idea of the proof is as follows. We separate the principal terms in operators $\mathbb{T},\,\mathbb{S}$. These terms are compact pseudodifferential operators of variable order. The singular values asymptotic for such operators is known, see \cite{K19}, \cite{K20}.  Then, using the asymptotic perturbation theory, we verify that remainder terms do not influence upon the obtained spectral asymptotics of principal terms.

\subsection{Operator $\mathbb{T}$ (mBM)}

Since a symbol of pseudodifferential operator should be smooth with respect to the dual variable $\xi$, we introduce an even smooth positive function ${\bf p}(\xi)$ such that
$$
{\bf p}(\xi) = |\xi|\quad \mbox{for}\quad |\xi|\ge2;\qquad {\bf p}(\xi)= 1\quad \mbox{for}\quad |\xi|\le 1.
$$

Now we consider the pseudodifferential operator of variable order 
 \begin{equation*}
 ({\mathbb{A}}f(x):=  C_*(H(x))\int\limits_{-\infty}^{\infty}e^{ix\xi}{\bf p}(\xi)^{-(H(x)+  \frac 12)}f(\xi)\,d\xi.
 \end{equation*}
 
 The singular value asymptotics for $\mathbb{A}$ is given by part 2 in Corollary \ref{Cor}, see Sec.5.1. Thus, formula (\ref{TT1}) is ensured by part 3 in Proposition \ref{p2} if we prove the following lemma. 
 
 \begin{Lemma}
 The following estimate holds:
	 \begin{equation}\label{t}
{\cal N}(t, {\mathbb{T}- \mathbb{A}}) \le C t^{\frac 1{\mathfrak{m}}-\mu},\quad t>1,
	\end{equation} 
with some $\mu>0$.
 \end{Lemma}

\textbf{Proof.} 
 The kernels of $\mathbb{T}$ and $\mathbb{A}$ have the same bounded multiplier $C_*(H(x))$, so we need to estimate singular values for the operator with the kernel
\begin{equation}\label{q}
\aligned
\frac{e^{ix\xi}-1}{|\xi|^{H(x)+\frac 12}}-\frac{e^{ix\xi}}{{\bf p}(\xi)^{H(x)+\frac 12}}= &\, q_1(x,\xi)-q_2(x,\xi)+q_3(x,\xi)\\
:= &\, \frac{\zeta(\xi)-1}{|\xi|^{H(x)+\frac 12}} -\zeta(\xi)\,\frac{e^{ix\xi}}{{\bf p}(\xi)^{H(x)+\frac 12}}+\zeta(\xi)\,\frac{e^{ix\xi}-1}{|\xi|^{H(x)+\frac 12}}.
\endaligned
\end{equation} 
Here $\zeta(\xi)$ is a fixed even cut-off function, 
$$
\zeta(\xi) = 0\quad \mbox{for}\quad |\xi| \ge3;\qquad \zeta(\xi) = 1\quad \mbox{for}\quad |\xi| \le2 
$$
(notice that ${\bf p}(\xi)=|\xi|$ for $\zeta(\xi)\ne 1$).

According to (\ref{q}), we need to estimate singular values of the operator ${\mathbb Q}_1-{\mathbb Q}_2+{\mathbb Q}_3$, where ${\mathbb Q}_j :L_2({\mathbb R}) \to L_2(0,1)$ are the integral operators with kernels $q_j(x,\xi)$, $j=1,2,3$. Due to the part 2 of Proposition \ref{p2} the estimate (\ref{t}) follows from similar estimates for operators ${\mathbb Q}_j$.

The estimate ${\cal N}(t,{\mathbb Q}_1)\le C(\varepsilon) \,t^{\frac 1{\mathfrak{m}+\lambda} +\varepsilon}$ for any $\varepsilon >0$ follows from Lemma \ref{L5}. 

The kernel $q_2(x,\xi)$ is bounded in $x$, smooth and compactly supported in $\xi$. Therefore, by Proposition \ref{p3}, ${\cal N}(t,{\mathbb Q}_2) \le C(\varepsilon) t^\varepsilon$ for any $\varepsilon>0$.

The kernel $q_3$ is singular at the point $\xi = 0$. We separate the principal part of this singularity
$$
q_3(x,\xi)= q_{3,0}(x,\xi)+ q_{3,1}(x,\xi),\qquad
q_{3,0}(x,\xi) :=\zeta(\xi)\,\frac {ix \xi}{|\xi|^{{H(x)+\frac 12}}}. 
$$  
So, we can write ${\mathbb Q}_3 =  {\mathbb Q}_{3,0}+{\mathbb Q}_{3,1}$. 

The function  $q_{3,1}(x,\cdot)$ is compactly supported and belongs to the Sobolev space $W_2^1(\mathbb{R})$ uniformly with respect to $x$. Therefore, Proposition \ref{p3} gives ${\cal N}(t,{\mathbb Q}_1)=O(t^{\frac 23})$. The required estimate follows from the inequality $\frac 1{\mathfrak{m}}>\frac 23$.

It remains to estimate singular values of the operator ${\mathbb Q}_{3,0}$.
Its kernel is nonzero only for $|\xi|<3$, and we can consider ${\mathbb Q}_{3,0}$ as the operator from $L_2(-3,3)$ to $L_2(0,1)$.

We introduce the isometry
$$
{\mathbb U}:\,L_2(-3,3)\to L_2(\mathbb R); \qquad ({\mathbb U}f)(z):=\sqrt{3}e^{-\frac {|z|}2}f(3\,{\rm sign}(z)e^{-|z|}).
$$
Then singular values of ${\mathbb Q}_{3,0}$ coincide with ones of ${\mathbb Q}_{3,0}{\mathbb U}^{-1}:L_2(0,\infty)\to L_2(0,1)$. Since the kernel $q_{3,0}$ is odd with respect to $\xi$, we have
\begin{equation}\label{Q+}
\big({\mathbb Q}_{3,0}{\mathbb U}^{-1}g\big)(x)=3^{1-H(x)} ix\int\limits_0^\infty\zeta(3e^{-z}) e^{-(1-H(x))z}\big(g(z)-g(-z)\big)\,dz.
\end{equation}
The function $1-H(\cdot)$ belongs to ${\cal C}^\beta[0,1]$ and is bounded away from zero. Therefore, the integrand in (\ref{Q+}) belongs to ${\cal C}^\beta$ in $x$ and decays exponentially as $z \to \infty$.  Lemma \ref{L5} yields the estimate ${\cal N}(t,{\mathbb Q}_{3,0}) \le C(\varepsilon) t^\varepsilon$ for any $\varepsilon>0$. 

Summing up the obtained estimates we arrive at (\ref{t}). 
\hfill $\square$

\subsection{Operator $\mathbb{S}$ (mfBM)}

First, we change the variable $z=y(1+w)$ in integral (\ref{K}) and rewrite the kernel of $\mathbb{S}$ as
\begin{equation}\label{K1}
K(x,y,H)= c_*(H)y^{H-\frac 12}\int\limits_{0}^{\frac {x-y}y}w^{H-\frac 32}(1+w)^{H-\frac 12}dw\,{\mathbb I}_{[0,x]}(y).
\end{equation}
Next, we separate the principal homogeneous term in (\ref{K1}) and write 
$$
K(x,y,H)=\widetilde q_1(x,y,H)+\widetilde q_2(x,y,H),
$$
where
\begin{eqnarray*}
\widetilde q_1(x,y,H) & := & c_*(H)y^{H-\frac 12}\int\limits_{0}^{\frac {x-y}y}w^{H-\frac 32}\,dw\,\chi_{[0,x]}(y)= \frac{c_*(H)}{H-\frac 12}\, (x-y)^{H-\frac 12}{\mathbb I}_{[0,x]}(y);\\
\widetilde q_2(x,y,H) & = & c_*(H)y^{H-\frac 12}\Phi\Big(\frac {x-y}y,H\Big){\mathbb I}_{[0,x]}(y),\\
\Phi(s,H) & := & \int\limits_{0}^s w^{H-\frac 32}\big((1+w)^{H-\frac 12}-1 \big)dw.
\end{eqnarray*}

Since $x,y \in (0,1)$, we can assume that $\widetilde q_1$ is multiplied by a cut-off function $\theta(x-y)$,  
\begin{equation}\label{theta}
\theta \in {\cal C}^\infty ({\mathbb R}),\qquad \theta(w)=1 \quad \mbox{for} \quad |w|\le 1, \qquad \theta (w)= 0 \quad \mbox{for}\quad |w|\ge 2.  
\end{equation}

Let ${\widetilde{\mathbb{Q}}_1}$ and ${\widetilde{\mathbb{Q}}_2}$ be operators in $L_2(0,1)$ with kernels $\widetilde q_1(x,y,H(x))$ and $\widetilde q_2(x,y,H(x))$ respectively. We claim that ${\widetilde{\mathbb{Q}}_1}$ is in fact a pseudodifferential operator of variable order. Indeed, we have
\begin{equation*}
(\widetilde{\mathbb{Q}}_1f)(x)= (2\pi)^{-1}c_*(H(x))\int\limits_{-\infty}^{\infty}\int\limits_0^1 e^{i(x-y)\xi}R\Big(\xi,H(x)-\frac 12\Big) f(y)\,dyd\xi,
\end{equation*}
where
 $$
 R(\xi,\gamma):=\int\limits_0^\infty e^{-iz\xi}\theta(z)\,\frac {z^\gamma}\gamma\, dz.
$$ 

For any $\gamma >-1$, $\gamma \ne 0$,  we have $R(\cdot,\gamma)\in{\cal C}^\infty(\mathbb{R})$. Moreover, up to a function of the Schwartz class ${\cal S}$, the function $R(\cdot,\gamma)$ coincides at infinity with the Fourier transform of $z_+^\gamma /\gamma$: 
$$
R(\xi,\gamma) = \Gamma (\gamma)|\xi|^{-(1+\gamma)}\exp\big(i\,\mbox{sign}(\xi)(1+\gamma) \pi/4 \big) +O(|\xi|^{-n}) \quad \mbox{for any}\quad n\in\mathbb{N}.
$$
Thus, $R(\xi,\gamma)$ is a classical symbol of order $-(1+\gamma)$. Therefore, $\widetilde{\mathbb{Q}}_1$ can be considered as a pseudodifferential operator of variable order, and the claim follows. 

We define a pseudodifferential operator 
\begin{equation*}
(\widetilde{\mathbb{A}}f)(x) := (2\pi)^{-1}c_*(H(x))\Gamma \big(H(x)-\frac 12\big)\int\limits_{-\infty}^{\infty}\int\limits_0^1 e^{i(x-y)\xi}\bigl( \widetilde{\bf p}(\xi)\bigr)^{-(\mathfrak{m}+h(x))}f(y)\,dyd\xi,
 \end{equation*} 
where $h(x)=H(x)-H_{min}$ while $\widetilde{\bf p}(\xi)$ is a smooth complex-valued function such that
$$
\widetilde{\bf p}(\xi)\ne0;\qquad
\widetilde{\bf p}(\xi)= |\xi|\exp\big(i\,\mbox{sign}(\xi) \pi/4\big)\quad\mbox{for}\quad |\xi|>1.
$$
The kernel of $\widetilde{\mathbb{Q}}_1 -\widetilde{\mathbb{A}}$ is smooth in $y$ and bounded in $x$. By Proposition \ref{p3} we have 
${\cal N}(t,\widetilde{\mathbb{Q}}_1-\widetilde{\mathbb{A}}) = O(t^{\varepsilon})$ as $t\to\infty$ for any $\varepsilon>0$, and part 3 of Proposition \ref{p2} gives  ${\cal N}(t,\widetilde{\mathbb{Q}}_1)={\cal N}(t,\widetilde{\mathbb{A}})+ O(t^\varepsilon)$.

Part 1 in Corollary \ref{Cor}, see Sec.5.1, gives the singular value asymptotics for the operator $\widetilde{\mathbb{A}}$ and therefore for the operator $\widetilde{\mathbb{Q}}_1$. Thus, formula (\ref{SS1}) is ensured by part 3 of Proposition \ref{p2} if we prove the following lemma. 

\begin{Lemma}
The following estimate holds:
\begin{equation}\label{s}
{\cal N}(t, \widetilde{\mathbb{Q}}_2) \le C t^{\frac 23},\quad t>1.
\end{equation}
\end{Lemma}

\textbf{Proof.} 
The kernel $\widetilde q_2$ has singularities on the diagonal $y=x$ and at the point $y=0$, i.e. for $s\equiv \frac {x-y}y=0$ and $s=\infty$ respectively. We consider the influence of these singularities upon the singular value asymptotic separately.

We consider two functions 
$$
\Phi_0(s,H)= \theta(s)\Phi(s,H),\qquad \Phi_1(s,H)= (1-\theta(s))\Phi(s,H),
$$
where the cut-off function $\theta$ is defined in (\ref{theta}), and denote by $\widetilde{\mathbb{Q}}_{2,j}$, $j=0,1$, the operators in $L_2(0,1)$ with kernels 
$$
\widetilde q_{2,j}(x,y,H(x))= c_*(H(x))y^{H(x)-\frac 12}\Phi_j\Big(\frac {x-y}y,H(x)\Big)\chi_{[0,x]}(y).
$$
Since $\widetilde{\mathbb{Q}}_2=\widetilde{\mathbb{Q}}_{2,0}+\widetilde{\mathbb{Q}}_{2,1}$, the estimate (\ref{s}) follows from similar estimates for operators $\widetilde{\mathbb{Q}}_{2,j}$.

We begin with the operator $\widetilde{\mathbb{Q}}_{2,0}$. The kernel $\widetilde q_{2,0}(x,y,H(x))$ does not vanish only for $\frac x3 < y< x$. Since $\Phi(s,H) = O(s^{H+\frac 12})$ and $\partial_s\Phi(s,H) = O(s^{H-\frac 12})$ as
$s \to 0$, we evidently have the estimate 
$$
|\widetilde q_{2,0}(x,y,H)| \le C\,\frac {(x-y)^{H+\frac 12}}y;\qquad |\partial_y\widetilde q_{2,0}(x,y,H)| \le C\,\Big(\frac {(x-y)^{H-\frac 12}}y+\frac {(x-y)^{H+\frac 12}}{y^2}\Big).
$$ 
We recall that the functional parameter satisfies $H(x)>\frac 12$. Therefore, for a fixed $x$ the function $\widetilde q_{2,0}(x,\cdot,H(x))$ is ${\cal C}^1$-smooth. Moreover, the following estimate holds:
$$\|\widetilde q_{2,0}(x,\cdot,H(x))\|^2_{W_2^1(0,1)} \le Cx^2\int\limits_{x/3}^x\Big(\frac{(x-y)^{H(x)-\frac 12}}{y^2}\Big)^2\,dy \le Cx^{2H(x)-2},
$$
so, the integral 
$$
\int\limits_0^1 \|\widetilde q_{2,0}(x,\cdot, H(x))\|_{W_2^1(0,1)}^2dx 
$$ 
converges. Now Proposition \ref{p3} yields the estimate 
$$
{\cal N}(t, \widetilde{\mathbb{Q}}_{2,0}) \le C t^{\frac 23},\quad t>1.
$$

Further, the kernel $\widetilde q_{2,1}(x,y,H(x))$ does not vanish only for $0\le y \le x/2$ and has singularity at the point $y=0$.  

Similarly to the estimate for the operator ${\mathbb Q}_{3,0}$ in the previous subsection, we introduce the isometry 
$$
\widetilde{\mathbb U}:\,L_2(0,1)\to L_2(0,\infty); \qquad (\widetilde{\mathbb U}f)(z):=e^{-\frac {z}2}f(e^{-z}).
$$
Then singular values of the operator $\widetilde{\mathbb Q}_{2,1}$ in $L_2(0,1)$ coincide with ones of the operator $\widetilde{\mathbb Q}_{2,1}\widetilde{\mathbb U}^{-1}:L_2(0,\infty)\to L_2(0,1)$. Changing the variable we obtain that the kernel of $\widetilde{\mathbb Q}_{2,1}\widetilde{\mathbb U}^{-1}$ is
$$
r(x,z,H(x)):= c_*(H(x))x^{H(x)} (1+s)^{-H(x)}\Phi_1(s)|_{s= xe^z-1}.
$$

The following estimates for $n\in\mathbb{N}\cup\{0\}$ are obvious:
\begin{equation*}
\big((1+s)\partial_s\big)^n\Phi_1(s,H)=O(s^{2H-1}),\qquad \mbox{as}\quad s \to \infty.
\end{equation*}
Since $\partial_z\big(g(xe^z-1)\big)=(1+s)\partial_sg(s)|_{s=xe^z-1}$, we obtain for any $n\in\mathbb{N}\cup\{0\}$
$$\big|\partial_z^n r(x,z,H(x))\big| \le C(n) x^{H(x)}(1+s)^{-(1-H(x))}|_{s=xe^z-1}\stackrel{*}{\le} C(n) e^{-(1-H(x))z}
$$ 
(the inequality $(*)$ follows from $H(x)>\frac 12$).

The function $1-H(\cdot)$ belongs to ${\cal C}^\beta[0,1]$ and is bounded away from zero. Therefore, the kernel $r(x,z,H(x))$ belongs to ${\cal C}^\beta$ in $x$ and decays exponentially as $z \to \infty$.  Lemma \ref{L5} yields the estimate ${\cal N}(t,\widetilde{\mathbb Q}_{2,1}) \le C(\varepsilon) t^\varepsilon$ for any $\varepsilon>0$. 

Summing up the estimates for $\widetilde{\mathbb{Q}}_{2,0}$ and $\widetilde{\mathbb{Q}}_{2,1}$ we arrive at (\ref{s}). 
\hfill $\square$

\section{Small ball asymptotics for mBM and mfBM}
\label{S:4}

As explained in the Introduction, one-term asymptotic of eigenvalues for covariance operator provides, under mild assumptions (see \cite[Theorem 1]{Naz09b}), the logarithmic $L_2$-small ball asymptotic for corresponding process. 
\medskip

We begin with the multifractional Brownian motion (\ref{mBM}). The case ${\rm meas}\, {\bf D}>0$ is in fact quite elementary. In this case formula (\ref{TT1}) and part 3 of Remark \ref{R:1} give
$$
{\cal N}(t,\mathbb{T})\sim\frac {\mathfrak{C}}{\pi}\, t^{\frac 1{\mathfrak {m}}},\quad \mbox{as}\quad t\to\infty,
$$
where
$$
\mathfrak{C}=
\big(\Gamma(2\mathfrak{m})\sin(\pi H_{min})\big)^{\frac 1{2\mathfrak{m}}}{\rm meas}\, {\bf D}.
$$
Since the function $t\mapsto {\cal N}(t,\mathbb{T})$ is in essense inverse to the function $k\mapsto s_k^{-1}(\mathbb{T})$, we have, as $k\to\infty$,
$$
s_k(\mathbb{T})\sim \Big(\frac {\mathfrak{C}}{\pi k}\Big)^{\mathfrak{m}} \quad\Longleftrightarrow \quad \lambda_k(\mathbb{T}\mathbb{T}^*)\sim \Big(\frac {\mathfrak{C}}{\pi k}\Big)^{2\mathfrak{m}}.
$$
Notice that $2\mathfrak{m}>1$. 
Applying Proposition 2.1 in \cite{NaNi04a} we obtain
$$
\lim_{\varepsilon \rightarrow 0}\varepsilon ^{\frac 2{2\mathfrak{m}-1}}\log \mathbb{P}\{||W^{H(\cdot)}||_{L_2(0,1)}\leq \varepsilon\}=-\ \frac {2\mathfrak{m}-1}{2}\left(\frac {\mathfrak{C}}{2\mathfrak{m} \sin(\frac
{\pi}{2\mathfrak{m}})}\right)^{\frac {2\mathfrak{m}}{2\mathfrak{m}-1}}.
$$

Now we consider the case ${\rm meas}\, {\bf D}=0$. In this case formula (\ref{TT1}) and part 3 of Remark \ref{R:1} give
$$
{\cal N}(t,\mathbb{T})\sim\frac {\widetilde{\mathfrak{C}}}{\pi}\, 
t^{\frac 1{\mathfrak{m}}}\,(\log t)^{-\sigma}  \, \varphi (\log t).
$$
where
$$
\widetilde{\mathfrak{C}}=
\big(\Gamma(2\mathfrak{m})\sin(\pi H_{min})\big)^{\frac 1{2\mathfrak{m}}}\mathfrak{m}^{2\sigma} \Gamma(\sigma+1).
$$
Therefore, we have, as $k\to\infty$,
$$
s_k(\mathbb{T})\sim \Big(\frac {\widetilde{\mathfrak{C}}}{\pi \mathfrak{m}^\sigma}\cdot\frac {\varphi (\log k)}{k\log^\sigma k}\Big)^{\mathfrak{m}} \quad\Longleftrightarrow \quad \lambda_k(\mathbb{T}\mathbb{T}^*)\sim \Big(\frac {\widetilde{\mathfrak{C}}}{\pi \mathfrak{m}^\sigma}\cdot\frac {\varphi (\log k)}{k\log^\sigma k}\Big)^{2\mathfrak{m}}.
$$

Since $\lambda_k$ is a sequence regularly varying with index $2\mathfrak{m}>1$, we can apply \cite[Theorem 4.2]{KNN} where 
a general situation was considered. Concretization of formula (4.5) in \cite{KNN} for our case gives
\begin{multline*}
\lim_{\varepsilon \rightarrow 0}\varepsilon ^{\frac 2{2\mathfrak{m}-1}}\bigg(\frac {\log^\sigma\frac 1{\varepsilon}}{\varphi(\log\frac 1{\varepsilon})}\bigg)^{\frac {2\mathfrak{m}}{2\mathfrak{m}-1}}\log \mathbb{P}\{||W^{H(\cdot)}||_{L_2(0,1)}\leq \varepsilon\}\\
=-\ \frac {2\mathfrak{m}-1}{2}\left(\frac {\widetilde{\mathfrak{C}}}{2\mathfrak{m} \sin(\frac
{\pi}{2\mathfrak{m}})}\Big(\frac {2\mathfrak{m}-1}{2\mathfrak{m}}\Big)^\sigma\right)^{\frac {2\mathfrak{m}}{2\mathfrak{m}-1}}\!\!.
\end{multline*}

Now we are able to formulate the final statement.

\begin{Theorem}
\label{T:2}
Assume that the variable Hurst parameter $0<H(x)<1$ satisfies the assumptions {\bf 1}--{\bf 3} before Theorem \ref{T:1} with $\tau$ subject to the condition (\ref{tau}). Then for the mBM $W^{H(\cdot)}$, the following relation holds as $\varepsilon\to0$: 
 \begin{enumerate}
\item If ${\rm meas}\, {\bf D}>0$ then
\begin{equation}\label{D>0}
\aligned
 & \log \mathbb{P}\{||W^{H(\cdot)}||_{L_2(0,1)}\leq \varepsilon\}\sim-\,\varepsilon ^{-\frac 1{H_{min}}}\\
\times &\, \frac {H_{min}\,{\rm meas}\, {\bf D}}{(2H_{min}+1)\sin\big(\frac
{\pi}{2H_{min}+1}\big)}\left(\frac {\Gamma(2H_{min}+1)\sin(\pi H_{min})\,{\rm meas}\, {\bf D}}{(2H_{min}+1)\sin\big(\frac
{\pi}{2H_{min}+1}\big)}\right)^{\frac 1{2H_{min}}}.
\endaligned
\end{equation}

\item If ${\rm meas}\, {\bf D}=0$ then
\begin{equation}\label{D=0}
\aligned
& \log \mathbb{P}\{||W^{H(\cdot)}||_{L_2(0,1)}\leq \varepsilon\}\sim-\, 
\varepsilon ^{-\frac 1{H_{min}}}
\bigg(\frac {\varphi(\log\frac 1{\varepsilon})}{\log^\sigma\frac 1{\varepsilon}}\bigg)^{\frac {2H_{min}+1}{2H_{min}}}\\
\times &\, H_{min}\left(\frac {\big(\Gamma(2H_{min}+1)\sin(\pi H_{min})\big)^{\frac 1{2H_{min}+1}}\Gamma(\sigma+1)\big(H_{min}(H_{min}+\frac 12)\big)^\sigma}{(2H_{min}+1)\sin\big(\frac {\pi}{2H_{min}+1}\big)} \right)^{\frac {2H_{min}+1}{2H_{min}}}\!\!.
\endaligned
\end{equation}
\end{enumerate}
\end{Theorem}

To illustrate this theorem, we give several examples. For simplicity only, we assume that $H_{min}=\frac 12$.
\medskip

{\bf Example 1}. Let $H(x)=\frac 12+(x-x_0)_+^\gamma$, $0<x_0\le 1$, $\gamma>0$. Then we have ${\bf D}=[0,x_0]$, and formula (\ref{D>0}) reads
$$\log \mathbb{P}\{||W^{H(\cdot)}||_{L_2(0,1)}\leq \varepsilon\}\sim-\,\frac{x_0^2}8\, \varepsilon^{-2}, \qquad \text{as}\quad \varepsilon\to0.
$$
For $x_0=1$ we obtain standard Wiener process on $[0,1]$, and this result is well known. However, for $x_0<1$ even this simplest result seems to be new.
\medskip

{\bf Example 2}. Let $H(x) = \frac 12 + |x-x_0|^\gamma$, $\gamma>0$. In this case we have purely power-like behavior of the measure of the small values set. Namely, formula (\ref{H0}) holds with
$$
\sigma = \frac 1{\gamma};\qquad \varphi (s)\equiv 2\quad \mbox{if}\quad 0<x_0<1; \qquad \varphi (s)\equiv 1\quad \mbox{if}\quad x_0=0,1.
$$
Therefore, formula (\ref{D=0}) gives
so 
\begin{equation}\label{ex1}
\log \mathbb{P}\{||W^{H(\cdot)}||_{L_2(0,1)}\leq \varepsilon\}\sim-\, \widehat C(x_0)\Gamma^2\Big(1+\frac 1\gamma\Big)\cdot \Big(\varepsilon \log^{\frac 1\gamma} \frac 1{\varepsilon} \Big)^{-2}, \qquad \text{as}\quad \varepsilon\to0,
\end{equation}
where
$$
\widehat C(x_0)=2^{-1-\frac 2{\gamma}}\quad \mbox{if}\quad 0<x_0<1; \qquad \widehat C(x_0)=2^{-3-\frac 2{\gamma}}\quad \mbox{if}\quad x_0=0,1.
$$
\medskip

{\bf Example 3}. Let $H(x)= \frac 12 +\min\limits_{1\le k \le N} |x-x_k|^{\gamma_k}$, with $0\le x_k\le 1$, $\gamma_k>0$, $k=1,\dots,N$. In this case the function $H$ attains minimum at several points but only the point(s) with maximal $\gamma_k$ affect the asymptotics. For instanse, if $\gamma_k=\gamma$ for $k\le n$, $\gamma_k<\gamma$ for $k> n$, and $0<x_k<1$ for $k\le n$ then formula (\ref{ex1}) holds with $\widehat C(x_0)=2^{-1-\frac 2{\gamma}}n^2$.
\medskip

{\bf Example 4}. Let $H(x)=\frac 12 +|x-x_0|^\gamma\log^b(|x-x_0|^{-1})$, with $x_0 \in (0,1)$, $\gamma >0$ and $b \in {\mathbb{R}}$. Then formula (\ref{H0}) holds with
$$
\sigma = \frac1\gamma; \qquad \varphi(s^{-1}) = 2\gamma^{\frac b\gamma } \big(\log s^{-1}\big)^{-\frac b\gamma} \big(1+o(1)\big) \qquad\text{as}\quad s\to0. 
$$
By part 4 of Remark \ref{R:1}, we can drop $o(1)$ in the latter relation, and formula (\ref{D=0}) gives
$$
\log \mathbb{P}\{||W^{H(\cdot)}||_{L_2(0,1)}\leq \varepsilon\}\sim-\,\frac{\gamma^{\frac {2b}\gamma}\Gamma^2(1+\frac 1\gamma)}{2^{1+\frac 2\gamma}}\cdot \Big(\varepsilon \log^{\frac 1\gamma} \frac 1{\varepsilon}\, \log^{\frac b\gamma}\log\frac 1{\varepsilon}\Big)^{-2}, \qquad \text{as}\quad \varepsilon\to0.
$$

{\bf Example 5}. Let $H(x)=\frac 12 +{\rm dist}^\gamma(x, \mathfrak{D})$, where $\mathfrak{D}$ is standard Cantor set. A tedious but simple calculation gives for small $s>0$
$$
{\rm meas}\,\{x\in [0,1] \mid 0<h(x)<s\} = s^{1-\frac{\log 2}{\log{3}}}\,\phi(\log s^{-1}),
$$
where $\phi$ is a {\it periodic} function. Therefore, in this case 
the assumption {\bf 3} before Theorem \ref{T:1} is not valid, and this case is not covered by our Theorem \ref{T:2}. We are planning to consider corresponding class of processes in a forthcoming paper.
\medskip

Now we turn to the multifractal Brownian motion (\ref{mfBM}). By part 1 of Remark \ref{R:1}, for $H_{min}>\frac 12$ formulae (\ref{TT1}) and (\ref{SS1}) coincide, and thus the logarithmic asymptotic for $X^{H(\cdot)}$ coincides with that of $W^{H(\cdot)}$.

\begin{Theorem}
\label{T:3}
Assume that the variable Hurst parameter $\frac 12<H(x)<1$ satisfies the assumptions {\bf 1}--{\bf 3} before Theorem \ref{T:1} with $\tau$ subject to the condition (\ref{tau}). Then for the mfBM $X^{H(\cdot)}$, the following relation holds as $\varepsilon\to0$: 
 \begin{enumerate}
\item If ${\rm meas}\, {\bf D}>0$ then
$$
\aligned
 & \log \mathbb{P}\{||X^{H(\cdot)}||_{L_2(0,1)}\leq \varepsilon\}\sim-\,\varepsilon ^{-\frac 1{H_{min}}}\\
\times &\, \frac {H_{min}\,{\rm meas}\, {\bf D}}{(2H_{min}+1)\sin\big(\frac
{\pi}{2H_{min}+1}\big)}\left(\frac {\Gamma(2H_{min}+1)\sin(\pi H_{min})\,{\rm meas}\, {\bf D}}{(2H_{min}+1)\sin\big(\frac
{\pi}{2H_{min}+1}\big)}\right)^{\frac 1{2H_{min}}}.
\endaligned
$$

\item If ${\rm meas}\, {\bf D}=0$ then
$$
\aligned
& \log \mathbb{P}\{||X^{H(\cdot)}||_{L_2(0,1)}\leq \varepsilon\}\sim-\, 
\varepsilon ^{-\frac 1{H_{min}}}
\bigg(\frac {\varphi(\log\frac 1{\varepsilon})}{\log^\sigma\frac 1{\varepsilon}}\bigg)^{\frac {2H_{min}+1}{2H_{min}}}\\
\times &\, H_{min}\left(\frac {\big(\Gamma(2H_{min}+1)\sin(\pi H_{min})\big)^{\frac 1{2H_{min}+1}}\Gamma(\sigma+1)\big(H_{min}(H_{min}+\frac 12)\big)^\sigma}{(2H_{min}+1)\sin\big(\frac {\pi}{2H_{min}+1}\big)} \right)^{\frac {2H_{min}+1}{2H_{min}}}\!\!.
\endaligned
$$
\end{enumerate}
\end{Theorem}

\begin{Remark}
For both processes (\ref{mBM}) and (\ref{mfBM}), in the case ${\rm meas}\, {\bf D} >0$ we obtain purely power asymptotics. Even in the case $\sigma<\nu$ when Theorem \ref{T:1} gives a two-term asymptotics, see part 3 in Remark \ref{R:1}, the estimate of the remainder term is not sufficient to obtain exact small ball asymptotics.

The asymptotic coefficient in this case depends on ${\rm meas}\, {\bf D}$ and for ${\rm meas}\, {\bf D}=1$ coincides with classical result of Bronski \cite{Bron}, see also \cite[Theorem 3.1]{NaNi04a}.
\end{Remark}

\section{Appendix}

\subsection{Asymptotics of singular values for pseudodifferential operators of variable order}

We consider a compact pseudodifferential operator ${\mathbb A}:L_2({\mathbb R})\to L_2(0,1)$,
\begin{equation}\label{AA} 
({\mathbb A}f)(x) := (2\pi)^{-1}a(h(x))\int\limits^{\infty}_{-\infty}\int\limits^{\infty}_{-\infty} e^{i(x-y)\xi}\bigl( {\bf p}(\xi)\bigr)^{-(m+h(x))}f(y)\,dyd\xi.
 \end{equation} 
Here  $m>\frac 12$, and $h\in{\cal C}^\lambda [0,1]$ is a non-negative function.  The complex-valued non-vanishing function ${\bf  p}$ is an elliptic symbol in H\"ormander class $S^1_{\rho,0}$, $0<\rho \le 1$, namely, 
\begin{eqnarray}
\label{p(xi)}
{\bf p}\in {\cal C}^\infty({\mathbb{R}}), && |{\bf p}(\xi)|= v|\xi|(1+O(|\xi|^{-\mu}))\quad\mbox{as}\quad |\xi|\to\infty \quad\mbox{for some}\quad \mu>0;
\\
\nonumber
&& \big| {\bf p}^{(n)} (\xi)\big| \le C(n)(1+|\xi)|)^{1-\rho\, n}\quad \mbox{for any} \quad n \in {\mathbb{N}}\cup 0.
\end{eqnarray}
The (complex-valued) multiplier $a\in{\cal C}^\infty[0,\infty)$ satisfies $a(0)\ne 0$.
\medskip
 
The order of decay of the symbol $p(\xi)^{-(m+h(x))}$ as  $|\xi| \to \infty$, which can be interpreted as the local order of operator ${\mathbb A}$, depends on $x$. So, we say that ${\mathbb A}$ is an operator of variable order.

Similarly to Section \ref{S:2}, we consider the set ${\bf D}:= \{x \in [0,1]\mid h(x)=0\}$.  

\begin{Theorem} \label{T:aux} 
Assume that the assumptions {\bf 1}--{\bf 3} before Theorem \ref{T:1} are fulfilled. Assume in addition that
\begin{equation*}
\tau>
\begin{cases}
\max\{0, \frac 12\,(1-\sigma)\}, & if \quad {\rm meas}\, {\bf D}>0;\\
\sigma m +\frac 12\,(1-\sigma), & if \quad {\rm meas}\, {\bf D}=0.
\end{cases}
\end{equation*}
Then, as $t \to \infty$,
\begin{equation}\label{TT1*}
{\cal N}(t,\mathbb{A})=\frac 1{\pi v}\,\big(a(0)\big)^{\frac 1m} \int\limits_0^1 t^{\frac 1{m+h(x)}}\,dx\,\big(1+O(\log^{-\nu}t)\big).
\end{equation}
Here $v$ is defined in (\ref{p(xi)}) while $\nu =\nu (m,\sigma,\tau)>0$. If ${\rm meas}\, {\bf D}>0$ then $\nu$ is arbitrary exponent less than $\frac{\tau -\frac 12\,(1-\sigma)}{m+1}$.
\end{Theorem}

\textbf{Proof.}  Assumptions {\bf 1}--{\bf 3} ensure that the conditions of \cite[theorem 1.3]{K20} are satisfied. In \cite{K20} the statement was proved, even in a multidimensional case, for $a\equiv 1$. The proof runs without changes for the multiplier $a(h_0(x))$. To include the nonsmooth multiplier $a(h(x))$ one needs to study the operator with the symbol 
$$
\bigl(a(h(x))-a(h_0(x))\bigr)\bigl( {\bf p}(\xi)\bigr)^{-(m+h(x))}.
$$
Repeating the argument of \cite[theorem 5.1]{K20} we obtain the estimate
$$
O(\log^{-\nu}t)\int\limits_0^1 t^{\frac 1{m+h(x)}}\,dx
$$
for the singular values counting function of this operator. Part 4 in Proposition \ref{p2} with regard to (\ref{Laplace}) completes the proof.
\hfill$\square$

\begin{Remark}\label{R:3}
\begin{enumerate}
\item Theorem \ref{T:aux} in fact claims that Weyl's spectral asymptotic formula 
$$
{\cal N}(t,\mathbb{A}) \sim  (2\pi)^{-1}meas\,\big\{ x,\xi \in (0,1) \times{\mathbb R} \mid  |a(h(x)){\bf p}(\xi)^{-(m+h(x))}| > t^{-1}\big\}
$$	
is valid for operator $\mathbb{A}$. Asymptotics (\ref{TT1*}) is just concretization of Weyl's formula for our problem with the remainder estimate.

\item If $h_1(x) \equiv 0$, then formula (\ref{TT1*}) is valid with $O(t^{-r})$ remainder estimate for some $r>0$.
\end{enumerate}
\end{Remark}

\begin{Corollary}\label{Cor}
\begin{enumerate}
 \item 
One can consider the operator in formula (\ref{AA}) as an operator acting in $L_2(0,1)$ if we extend $u\in L_2(0,1)$ by zero. Formula (\ref{TT1*}) remains valid in this case.

\item

Since Fourier transform is the unitary operator in $L_2({\mathbb R} )$, formula (\ref{TT1*}) holds also for the operator 
\begin{equation*}
({\mathbb{A}}f)(x) := (2\pi)^{-\frac 12}a(h(x)) \int\limits^\infty_{-\infty} e^{ix\xi}\bigl( {\bf p}(\xi)\bigr)^{-(m+h(x))}f(\xi)d\xi.
\end{equation*}

\end{enumerate}

\end{Corollary}

\subsection{Estimates for singular values of integral operator}

We need the following results in asymptotic perturbations theory.  

\begin{Proposition}\label{p2}
\begin{enumerate}
\item
Let $\mathbb{A}:{\cal H}_1 \to {\cal H}_2$, $\mathbb{B}: {\cal H}_2 \to {\cal H}_3$ be compact operators. If ${\cal N}(t,\mathbb{A}) = O(t^{p_1})$ and  

${\cal N}(t,\mathbb{B}) = O(t^{p_2})$ as $t \to \infty$ then ${\cal N}(t,\mathbb{BA}) = O(t^{p})$ with $\frac 1p=\frac 1{p_1}+ \frac 1{p_2}$.

\item
Let $\mathbb{A},\,\mathbb{B}: {\cal H}_1 \to {\cal H}_2 $ be compact operators. If ${\cal N}(t,\mathbb{A}) = O(t^p)$ and  ${\cal N}(t,\mathbb{B}) = O(t^p)$ as $t \to \infty$ then ${\cal N}(t,\mathbb{A+B}) = O(t^p)$.

\item 
Let $\mathbb{A},\,\mathbb{B}: {\cal H}_1 \to {\cal H}_2 $ be compact operators. If ${\cal N}(t,\mathbb{A}) = t^pV(t)$ as $t \to \infty$ with a slowly varying function $V$, and ${\cal N}(t,\mathbb{B}) = O(t^q)$ for some $q<p$, then
$$
{\cal N}(t,\mathbb{A+B}) = t^p \,(V(t)+ O(t^{-r})),\quad \mbox{where}\quad  r<\frac {p-q}{q+1}.
$$
\item
Let $\mathbb{A},\,\mathbb{B}: {\cal H}_1 \to {\cal H}_2 $ be compact operators. Assume that 

$$
{\cal N}(t,\mathbb{A}) = t^p(\delta + \varkappa \log^{-a}t\,V(\log t)) \quad\mbox{as}\quad t\to\infty
$$
with $a>0$ and a slowly varying function $V$, and ${\cal N}(t,\mathbb{B}) = O(t^p\log^{-b} t)$ for some $b>a$. Then 
$${\cal N}(t,\mathbb{A}+\mathbb{B})= t^p(\delta +\varkappa \log^{-a}t\,V(\log t)+ O(\log^{-r}t)),
$$

where $r$ is arbitrary exponent such that
$$
r<\frac b{p+1},\quad\mbox{if}\quad \delta>0;\qquad
r<\frac {ap+b}{p+1},\quad\mbox{if}\quad \delta=0.
$$
\end{enumerate}

\end{Proposition}
The assertions 1 and 2 are elementary consequences of the inequalities (17) and (19) in \cite[Sec. 11.1]{BS10}, respectively. The statements 3 and 4 refine the estimate in \cite[Lemma 3.1]{BS01} and are contained in \cite[Lemma 2.1]{K20}. 
\medskip

Now we formulate the singular values estimates for integral operator in terms of the properties of its kernel. The results are given for $d$-dimensional domains, though we need only the case $d=1$. In what follows $W_2^{\lambda}$ stands for the standard Sobolev--Slobodetskii space, see, e.g., \cite[Sec. 2.3.1]{Tr}.

Let $\Omega$ be a bounded domain in ${\mathbb R}^d$ with Lipschitz boundary. We put
 \begin{equation*}
 \mathfrak{R}:L_2({\mathbb R}^d)\to L_2(\Omega),\qquad
 (\mathfrak{R}f)(x):= a(x)\int\limits_{{\mathbb R}^d}R(x,y)f(y)\,dy, 
 \end{equation*}
where $a\in L_\infty(\Omega)$. 
 
\begin{Proposition}\label{p3} (see \cite[\S 11.8]{BS10}). 
Let the function $R(\cdot,y)\in W_2^{\lambda}(\Omega)$ for a.e. $y \in {\mathbb R}^d$, and let
\begin{equation}\label{ss3}
M^2 :=\int\limits_{{\mathbb R}^d} \|R(\cdot, y)\|^2_{W_2^\lambda(\Omega)}dy  < \infty.
\end{equation}	
Then the following estimate  holds:
\begin{equation*}
{\cal N}(t,\mathfrak{R})\le C(\Omega) \big(M \|a\|_{L_\infty (\Omega)}\big)^pt^p, \qquad \frac 1p = \frac{1}{2} + \frac{\lambda}{d}.
\end{equation*}
\end{Proposition}

Now we take into account the decay of the kernel with respect to the second variable. Let
$$
\widetilde{\mathfrak{R}}: L_2({\mathbb R}^d)\to L_2(\Omega),\qquad (\widetilde{\mathfrak{R}}f)(x):= \int\limits_{{\mathbb R}^d}\widetilde R(x,\xi)f(\xi)\,d\xi.
$$
\begin{Lemma}\label{L5}
Let the function $\widetilde R(\cdot,\xi)\in W_2^\lambda(\Omega)$ for a.e. $\xi \in {\mathbb R}^d$, and let
\begin{equation}\label{r1}
\|\widetilde R(\cdot, \xi)\|_{W_2^\lambda(\Omega)} \le M(\ell,\lambda)(1+|\xi|)^{-\ell} \qquad \mbox{for some}\quad \ell > \frac d2.
\end{equation}
Assume also that $\widetilde R(x,\cdot)\in {\cal C}^\infty(\mathbb{R}^d)$, and
\begin{equation*}
|\partial^\alpha_\xi \widetilde R(x,\xi)| \le C(\alpha)(1+|\xi|)^{-(\ell+ \rho|\alpha|) },\qquad x\in \Omega,\ \xi\in {\mathbb R}^d
\end{equation*}
for some $\rho >0$ and any multi-index $\alpha$.

Then
$$
{\cal N}(t,\widetilde{\mathfrak{R}}) \le C(\ell,\lambda,p,\Omega)t^p \qquad \mbox{for any}\quad p> \frac d{\ell+\lambda}.
$$
\end{Lemma}

\textbf{Proof.} Since the Fourier transform $\mathfrak{F}:L_2({\mathbb R}^d) \to L_2({\mathbb R}^d)$ is unitary, the singular values of operators $\widetilde{\mathfrak{R}}$ and $\mathfrak{R}= \widetilde{\mathfrak{R}}\,\mathfrak{F}$ coincide. 

We split $\mathfrak{R}$ into two parts:
$$
\mathfrak{R}=\mathfrak{R}_0 +\mathfrak{R}_1,
$$
with the kernels 
$$
R_0(x,y) = \theta(|y|)R(x,y),\qquad R_1(x,y)=(1-\theta(|y|))R(x,y),
$$
where the cut-off function $\theta$ is defined in (\ref{theta}), while
$$R(x,y) = (2\pi)^{-\frac d2}\int\limits_{{\mathbb R}^d}\widetilde R(x,\xi)e^{-iy\xi}d\xi, \quad x\in \Omega,\ y\in {\mathbb R}^d.
$$ 
Integrating by parts with respect to $\xi$ and using the identity
$$
e^{-iy\xi}= (-\Delta_\xi)^N|y|^{-2N}e^{-iy\xi},
$$
we obtain that the function $R_1(x,\cdot)$ belongs to the Schwartz class uniformly with respect to $x \in \Omega$. 
By \cite[theorem 4.8]{BS77}, this gives ${\cal N}(t,\mathfrak{R}_1)=O(t^\varepsilon)$ for any $\varepsilon>0$. So, part 2 of Proposition \ref{p2} shows that the estimate of ${\cal N}(t,\mathfrak{R})$ is governed by the estimate of ${\cal N}(t,\mathfrak{R}_0)$. 

We choose some $0<l <\ell-\frac d2$ and
write $\mathfrak{R}_0$ in the form $\mathfrak{R}_0=\mathfrak{R}_{0,0}\mathfrak{R}_{0,1}$, where the kernels of $\mathfrak{R}_{0,j}$, $j=0,1$, are
\begin{equation*}
\aligned
R_{0,0}(x,\xi)= &\ (2\pi)^{-\frac d2}\widetilde R(x,\xi)(1+|\xi|^2)^{\frac l2}, && x\in \Omega,\ \xi\in {\mathbb R}^d;\\
R_{0,1}(\xi,y)= &\ e^{-iy\xi}(1+|\xi|^2)^{-\frac l2}\theta(|y|), && \xi,y\in {\mathbb R}^d.
\endaligned
\end{equation*} 
For operator  $\mathfrak{R}_{0,0}$ 
we have the estimate ${\cal N}(t,\mathfrak{R}_{0,0})\le C(\Omega)t^{\frac dl}$ (this is a particular case of the Rozenblum--Lieb--Cwikel estimate, see, e.g., \cite[theorem 6.5]{BS77}). The estimate (\ref{r1}) ensures the inequality (\ref{ss3}), and Proposition \ref{p3} yields the estimate ${\cal N}(t,\mathfrak{R}_{0,0})\le C(\ell,\lambda,l)t^{p_1}$, $\frac 1{p_1}=\frac 12 +\frac {\lambda}d$. By part 1 of Proposition \ref{p2}, we have
$$
{\cal N}(t,\mathfrak{R}_0)\le C(\ell,\lambda,l,\Omega) t^p,\qquad \frac 1p=\frac 12 + \frac {\lambda}d + \frac ld < \frac {\ell+\lambda}d,
$$  
and the statement follows.
\hfill $\square$

\end{document}